\documentclass{article}
\usepackage{amsmath,amsthm,amssymb}

\newcommand{\gd}{\delta}
\newcommand{\gw}{\omega}

\newcommand{\gs}{\sigma}

\newcommand{\R}{\mathbb{R}}

\newcommand{\meager}{\mathtt{meager}}
\newcommand{\lnull}{\mathtt{null}}

\newtheorem{theorem}{Theorem}[section]
\newtheorem{lemma}[theorem]{Lemma}
\newtheorem{claim}[theorem]{Claim}
\newtheorem{corollary}[theorem]{Corollary}
\newtheorem{fact}[theorem]{Fact}

\theoremstyle{definition}

\newtheorem{question}[theorem]{Question}
\newcommand{\B}{{\mathcal B}}

\newcommand{\F}{{\mathcal F}}
\newcommand{\N}{{\mathbb N}}

\newcommand{\e}{\varepsilon}

\newcommand{\bbP}{\mathbb P}

\newcommand{\co}{\complement}

\newcounter{my_enumerate_counter}

\newcommand{\lbl}{\label}

\DeclareMathOperator{\Null}{Null}
\title{Between Maharam's and von Neumann's 
problems}
\author{
Ilijas Farah
\thanks{Partially supported by NSERC.}\\
York University\\
\and
Jind{\v r}ich Zapletal
\thanks {Partially supported by grant GA {\v C}R 
201-03-0933, NSF grant DMS 0071437, and  the first author's NSERC grant.
Results of this paper were obtained in November and December 2003
while the second author was visiting 
York University.}\\
University of Florida}

\begin{document}
\maketitle

\begin{abstract}
We show  that in the definable context the Maharam and von Neumann problems essentially
coincide.
We also prove that the random forcing is the only definable ccc forcing adding a single real that 
does not make the ground model reals null, and that the only pairs of definable ccc $\sigma$-ideals
with  the Fubini property are $(\meager,\meager)$ and $(\lnull,\lnull)$. 
\end{abstract}

%\section{The theorem}

In  Scottish Book, von Neumann asked whether every ccc, weakly distributive 
complete Boolean algebra carries a strictly positive probability measure. 
In her commentary to von Neumann's problem (\cite{Mah:vonNeumann}) D. Maharam pointed out
that this problem naturally splits into two problems: 
(a) whether all such algebras carry a strictly positive continuous submeasure, 
and (b) whether every algebra that carries a strictly positive continuous 
submeasure carries a strictly positive measure. 
The latter problem is known under the names of Maharam's Problem and Control Measure 
Problem (see \cite{Ka:Maharam}, \cite[\S 393]{Fr:MT3}). 
While von Neumann's problem has a consistently negative answer (\cite{Mah:Algebraic}), 
Maharam's problem can be stated as a $\mathbf \Sigma^1_2$ statement and is therefore by Shoenfield's 
theorem
absolute between transitive models of set theory containing all countable ordinals. 

Very recently
Balcar, Jech and Paz\'ak  announced  (\cite{BJP}) 
that under the   P-ideal dichotomy (\cite{To:Dichotomy}) 
 every c.c.c. weakly distributive complete Boolean 
 algebra carries a strictly positive Maharam submeasure. 
About a month before their
announcement, we discovered a similar result for suitably definable forcings
(see below for definitions).

\begin{theorem}\label{maintheorem}
Let $I$ be a  c.c.c. $\gs$-ideal on Borel subsets of $2^\gw$ that is 
analytic on $G_\delta$. The following are equivalent:

\begin{itemize}
\item $P_I$ is a weakly distributive notion of forcing
\item there is a Maharam submeasure on $2^\gw$ such that $I$ is the $\gs$-ideal of its null sets.
\end{itemize}
A suitable large cardinal assumption implies 
that the assumption that $I$ is analytic on $G_\delta$ can be relaxed to 
`$I$ is  definable.'
\end{theorem}

By a result of Shelah (\cite{Sh:480}), the first clause is equivalent to saying that~$P_I$ 
does not add Cohen reals, so we have a dichotomy for definable ccc 
forcing notions of the 
form $P_I$. 
Here $P_I$ is the partial ordering of $I$-positive Borel sets under inclusion. 
Thus, the von Neumann's problem restricted to definable partial orders of the form $P_I$
coincides with the Control Measure Problem. By an absoluteness argument, 
this also follows from Balcar--Jech--Paz\'ak's result. 

We were able to draw an  interesting consequence of the theorems. 
In order to state it succintly we quote the large cardinal version.

\begin{corollary}[LC] \label{C.random.0}
Suppose  $I$ is a  definable c.c.c. $\gs$-ideal on $2^\gw$. 
Then exactly one of the following holds: 
\begin{enumerate}
\item There is a Borel set $B\subset \R\times\R$
 with all  vertical sections in $I$ and all horizontal sections of 
 full Lebesgue measure. 
\item There is a condition $p\in P_I$ such that $P_I$ below $p$ 
is isomorphic to the random forcing.
\end{enumerate}
\end{corollary}

In other words, if $P_I$ does not force that the set of ground model reals is null, 
then $P_I$ is the random forcing.  Modulo Theorem~\ref{maintheorem}, this is really 
a consequence of a result of Christensen (\cite{Chr:Fubini}). 
There is a curious duality with an earlier result of Shelah 
that a similar result holds on the meager side. 

\begin{fact} [LC] (\cite{Sh:630}, see also \cite{z:products})
Suppose  $I$ is a  definable c.c.c. $\gs$-ideal on $2^\gw$. 
Then exactly one of the following holds: 
\begin{enumerate}
\item There is a Borel set $B\subset \R\times\R$
 with all  vertical sections in $I$ and all horizontal sections comeager. 
\item There is a condition $p\in P_I$ such that $P_I$ below $p$ 
is isomorphic to the Cohen forcing.
\end{enumerate}
\end{fact}

Another  attractive corollary is that the only
definable c.c.c. $\gs$-ideals for which Fubini theorem holds 
are the meager and null ideal (Theorem~\ref{T.non-commuting}). 
This shows that the only `reasonable' ideals as introduced by Kunen in \cite{Ku:Random-Cohen}
are  $\meager$ and $\lnull$.

\subsubsection*{Terminology} The notation in the paper
 follows the set theoretic standard of \cite{jech:set}. 
An ideal $I$ is analytic on $G_\delta$ if for every $G_\delta$ set 
$A\subseteq 2^\omega\times 2^\omega$ the set of all $x$ such that the 
vertical section of $A$ at $x$ is in $I$ is analytic.
Note that both $\meager$ and $\lnull$ are analytic on $G_\delta$ (see \cite{Ke:Classical}). 
Throughout the paper we will say that an ideal is  \emph{definable} if 
it belongs to $L(\R)$. 
The suitable large cardinal assumption in Theorem~\ref{maintheorem} if $I$ is in $L(\R)$ 
is that there are $\omega$ Woodin cardinals with a measurable above them all. 
In all the subsequent results of this note no large cardinal assumptions are needed if
$I$ is assumed to be analytic on $G_\delta$. 
 For large cardinals and $L(\R)$ see e.g., \cite{Kana:Book}.

A \emph{Maharam submeasure} (or a \emph{continuous submeasure})
on a complete Boolean algebra $\B$ is a function 
$\phi$ such that 
\begin{enumerate}
\item   $A\subseteq B$ implies $\phi(A)\leq \phi(B)$, 
\item $\phi(A\cup B)\leq \phi(A)+\phi(B)$, 
\item $\phi(0_\B)=0$, and
\item if $A_n$ is a decreasing sequence in $\B$ then 
$\phi(\bigcap_n A_n)=\lim_n \phi(A_n)$. 
\end{enumerate}
An algebra that carries a strictly positive Maharam submeasure is called a \emph{submeasure algebra}. 

A forcing notion is \emph{bounding} (or \emph{weakly distributive}) if every element of  
 $\omega^\omega$ in the extension 
is dominated by a ground-model function in $\omega^\omega$.
We use the words ``bounding'' and ``weakly distributive'' interchangeably. It is entirely irrelevant
which uncountable Polish space the ideals in question measure; our choice is 
the Cantor space $2^\gw$ for definiteness and ease of notation. To weed out
trivial cases, we assume that ideals contain all singletons.

\section{The proof of Theorem~\ref{maintheorem}}

It is a standard fact that if $I$ is the ideal of null 
sets of some Maharam submeasure, then the poset
$P_I$ is weakly distributive (see e.g., \cite[392I]{Fr:MT3}).
Suppose now that $I$ is a definable, weakly distributive 
c.c.c. $\gs$-ideal on Borel subsets of~$2^\gw$. 
To find a Maharam submeasure generating the $\gs$-ideal $I$ 
we will use two ingredients. One is almost trivial:

\begin{fact}[\cite{z:book} Lemma 2.2.3.]
\label{zfact}
Suppose that $I$ is a $\gs$-ideal on $2^\gw$ such that $P_I$ is 
proper. The following are equivalent:

\begin{itemize}
\item $P_I$ is weakly distributive
\item compact sets are dense (every $I$-positive Borel set has an $I$-positive compact subset)
and continuous reading of names (for every $I$-positive 
Borel set $B$ and a Borel function $f:B\to\gw^\gw$
there is an $I$-positive set $C\subset B$ such that $f\restriction C$ is continuous).
\end{itemize}
\end{fact}

\noindent Note that this implies that our ideal $I$ has a basis consisting of $G_\gd$ sets. For let
$A\in I$ be a Borel set. The collection of compact
$I$-positive sets disjoint from $A$ is dense in $P_I$: for every $I$-positive Borel set $B$, 
the set $B\setminus A$ is still Borel and $I$-positive and therefore it 
has a compact $I$-positive subset.
Choose then a maximal antichain $X$ consisting of such compact sets. Since 
the poset $P_I$ is c.c.c.,
it is the case that the antichain $X$ is countable, and $2^\gw\setminus \bigcup X$ is the required
$G_\gd$ set in the ideal $I$ covering the set $A$.

The other ingredient is a result of Solecki. 
For an ideal $I$ on $2^\omega$ let  $\hat I$ be the collection
of subsets of $2^{<\gw}$ defined by putting $a\in\hat I$ if the set
$B_a=\{r\in 2^\gw:$ for infinitely many $n, r\restriction n\in a\}$
is in  $I$.
It is immediate that $\hat I$ is an 
ideal, because $B_{a\cup b}=B_a\cup B_b$ and so if both $B_a, B_b\subset 2^\gw$
are in the $\gs$-ideal $I$, so is $B_{a\cup b}$.

\begin{fact}
\label{sfact}
 Suppose that $I$ is a  $\gs$-ideal on Borel subsets of $2^\gw$ that is 
 analytic on $G_\delta$. The following are equivalent:
\begin{itemize}
\item $\hat I$ is a $P$-ideal and $I$ has a basis consisting of $G_\gd$ sets
\item There is a Maharam submeasure on $2^\gw$ such that $I$ is the collection of its null sets.
\end{itemize}
Furthermore, large cardinals imply that this equivalence holds for every definable ideal $I$.  
\end{fact}
 
\begin{proof} This was proved in \cite[Theorem 5.2]{Sol:AnalyticII}.
in the case when $I$ is 
analytic on~$G_\delta$. The definability assumption was used in this proof only to show that 
$\hat I$ is analytic, and 
assuming large cardinals, in \cite[Theorem~4]{To:Definable}
it was proved that all definable P-ideals are analytic. 
\end{proof} 

 Fact~\ref{sfact} clearly implies that we will be done once we prove that
$\hat I$ is a $P$-ideal.
 To verify that this holds, fix a collection
$\{a_n:n\in\gw\}\subset\hat I$ and aim to construct some set $b\in\hat I$ which 
includes each of them up to
a finite set. 

\begin{claim}
The collection of compact $I$-positive sets $C$ such that their associated 
tree on $2^{<\gw}$ has a finite intersection
with each set $a_n:n\in\gw$, is dense in~$P_I$.
\end{claim}

\begin{proof}
Suppose $A\in P_I$ is a positive Borel set, and let 
$B=A\setminus\bigcup_n B_{a_n}$. $B$ is still an $I$-positive
Borel set, and the function $f:B\to\gw^\gw$, 
$f(r)(n)=\max\{m\in\gw:r\restriction m\in a_n\}$, is Borel
and well-defined on it. By Fact~\ref{zfact}, there is 
an $I$-positive compact set $C\subset B$ such that
$f\restriction C$ is continuous. By a compactness argument, 
for every number $n$ the set $\{f(r)(n):r\in C\}$
is finite. The claim follows.
\end{proof}

Let $X$ be a maximal antichain of $I$-positive compact 
sets from the claim. Since the poset $P_I$ is c.c.c.,
the collection $X$ is countable, enumerated as 
$\{C_k:k\in\gw\}$. Let $T_k\subset 2^{<\gw}$ be the tree
associated with the compact set $C_k$ for every $k\in\gw$. 
Finally, let $b\subset 2^{<\gw}$ be the set
$\bigcup_n (a_n\setminus\bigcup_{k\in n}T_k)$. It is clear 
that the set $b$ includes every $a_n$ up to
a finite piece. To show that $B_b\in I$ and $b\in\hat I$, 
note that for every number $k\in\gw$ the intersection
$T_k\cap b=\bigcup_{n\leq k} a_n$ is finite, and so the 
set $B_b$ is disjoint from $\bigcup X$. However,
the antichain $X\subset P_I$ was chosen to be maximal, 
and therefore the set $2^\gw\setminus\bigcup X$ is $I$-small
and so is its subset $B_b$. The theorem follows.

\section{Fubini failing}

A submeasure $\phi$ is \emph{pathological} if it does not dominate 
a positive nonzero finitely additive
functional. A \emph{control measure} for a Maharam submeasure $\phi$ is a measure $\mu$
that has the same null sets as $\phi$. A Maharam submeasure is \emph{Borel} if it is defined
on the Borel algebra on $2^\omega$. 
A submeasure is \emph{diffuse} if all countable sets are null. 
All results of this section are probably well-known.

\begin{lemma} \lbl{L.path} The following are equivalent for a diffuse Borel 
Maharam submeasure $\phi$. 
\begin{enumerate}
\item\lbl{L.path.1} $\phi$ is pathological. 
\item\lbl{L.path.2} There is a $\phi$-positive set $B$ such that the restriction of $\phi$ to 
$B$ has a control measure. 
\item\lbl{L.path.3} There is a $\phi$-positive set $B$ such 
that $P_{\Null(\phi)}$ is forcing equivalent to 
random below $B$. 
\end{enumerate}
\end{lemma} 

\begin{proof} Let us write $I=\Null(\phi)$. 
Assume \eqref{L.path.1}, so 
there is a nonzero finitely additive functional $\nu\leq\phi$ 
dominated by $\phi$. There are two cases. 

Assume  there is a $\phi$-positive set $B$ such that $\nu(C)\neq 0$ for 
every $I$-positive  set $C\subset B$. Then Borel$/I$ is 
weakly distributive (see e.g., \cite[392I]{Fr:MT3}). 
By 
 \cite[391D]{Fr:MT3}
there is a strictly positive measure on Borel$/I$, and therefore \eqref{L.path.2} holds. 
 
Otherwise, 
every $\phi$-positive  set  $B$ contains a $\phi$-positive  set  $C$ such 
that $\nu(C)=0$. In this case, choose
a maximal antichain $\{C_n:n\in\gw\}$ of sets such that $\phi(C_n)>0$ and $\nu(C_n)=0$, 
 enumerated using the 
ccc of Borel$/I$. Consider
the sets $D_m=\bigcup_{n>m}C_n$. By the finite additivity of the functional $\nu$ it is the
 case that
$\nu(D_m)=\nu(2^\gw)$ for all $m\in\gw$. By the continuity of the submeasure $\phi$, the numbers
$\phi(D_m)$ converge to zero, since the $D_m$'s for a decreasing collection of sets with 
empty intersection.
This however contradicts the assumption $\nu\leq\phi$.

Clause \eqref{L.path.2} implies \eqref{L.path.1} by \cite[Theorem~2]{HeChr:Example}. 
The equivalence of \eqref{L.path.2} and \eqref{L.path.3} follows by Maharam's theorem. 
\end{proof}

\begin{lemma} \lbl{L.decompose}
If $\phi$ is a Borel Maharam submeasure 
then there is a Borel set  $A$ such that 
$\phi$ has a control measure on $B$ and is pathological on $B^\co$. 
\end{lemma} 

\begin{proof} Find a maximal family $\F$ of pairwise orthogonal measures dominated by $\phi$, 
and let $B$ be the union of their supports. By the ccc-ness of 
Borel/$\Null(\phi)$, $\F$ is countable. 
If $\F=\{\mu_i| i\in \omega\}$ then 
$\sum_i 2^{-i}\mu_i$ is a control measure for $\phi$ on $B$. 
By Lemma~\ref{L.path}, $\phi$ is pathological on the complement of $B$. 
\end{proof} 

Lemma~\ref{L.Chr} below was roughly proved by  Christensen \cite[Theorem~6]{Chr:Fubini}.
 We shall use his result. 
Let $\mu$ denote the Lebesgue measure on $[0,1]$; the choice is immaterial as any other diffuse
Borel probability measure would do. 

\begin{lemma}\label{L.Chr}
 Suppose  $I$ is  the null ideal for some  Maharam Borel submeasure $\phi$ on $2^\omega$. 
Exactly one of the following holds:
\begin{enumerate}
\item \lbl{L.Chr.1}There is a $\phi$-positive Borel set $B$ such that the restriction of $\phi$ to 
$B$ has a control measure. 
\item\lbl{L.Chr.2}  There is a Borel set $C\subseteq [0,1]\times 2^\omega$ such that 
$\phi(C_x)=0$ for all $x\in [0,1]$ and $\mu([0,1]\setminus C^y)=0$ for every $y\in 2^\omega$. 
\end{enumerate}
\end{lemma}

\begin{proof} By Fubini's theorem,  \eqref{L.Chr.1} excludes \eqref{L.Chr.2}. 
 Suppose now that \eqref{L.Chr.1}  fails. 
By Lemma~\ref{L.path}, $\phi$ is pathological. 
Christensen proved in \cite[Theorem 6]{Chr:Fubini}, Theorem 6 
that if $\phi$ is pathological
then \eqref{L.Chr.2} holds. 
\end{proof}

A submeasure $\phi$ on $2^\omega$ is \emph{normalized} if $\phi(2^\omega)=1$. 

\begin{lemma} \lbl{L.KRL}
Assume $\psi$ is a normalized pathological Borel submeasure. 
Then for every $n\in \N$ there are pairwise disjoint sets $A_i$ ($i<n$) 
of submeasure at least $1/3$ each.  
\end{lemma} 

\begin{proof} This was proved by Kalton and Roberts (\cite{KaRo:Uniformly}) for an unspecified 
 $\e>0$ in place of $1/3$, and sharpened by  Louveau  (\cite{Lo:Progres}) to the present form. 
 \end{proof} 

 \begin{lemma} \label{L.non-Fubini}
 Assume $\phi$ and $\psi$ are normalized diffuse Borel 
 Maharam submeasures on $2^\omega$ and $\psi$ 
is pathological.  
  Then there is a Borel set $C\subseteq 2^\omega\times 2^\omega$
 such that $\psi(C_x)\geq 1/3$ for all $x\in 2^\omega$ and $\phi(C^y)=0$ for all 
 $y\in 2^\omega$. 
 \end{lemma} 
 
 \begin{proof} Since $\phi$ is diffuse and Maharam, every set of submeasure $\delta$ has a subset of 
 submeasure $\epsilon$ for every $\epsilon\in [0,\delta]$. For each $n$ a maximal antichain 
 of Borel sets such that the submeasure of each one is between $2^{-n-1}$ and $2^{-n}$. 
 Since $\phi$ is  Maharam, this antichain is finite and we can enumerate it as 
 $B^n_i$ ($i<k_n$). 
 Using  Lemma~\ref{L.KRL}, 
  fix a partition of $2^\omega$ into Borel  sets $A^n_i$ ($i<k_n$) 
 such that $\psi(A^n_i)\geq 1/3$ for all $n$ and $i$. 
 Let 
 \begin{align*}
 C(n)&=\bigcup_{i=0}^{k_n-1} B^n_i\times A^n_i\\
 C&=\bigcap_{m=0}^\infty \bigcup_{n=m}^\infty C(n). 
 \end{align*}
Note that $\psi(C(n)_x)\geq 1/3$ and that $\phi(C(n)^y)\leq 2^{-n}$ for all $x,y$ in $2^\omega$. 
Therefore for all $x,y$ we have 
$\psi(C_x)\geq 1/3$ and $\phi(C^y)\leq \sum_{n=m}^\infty 2^{-n}=2^{-m+1}$ for all $m$, hence
$\phi(C^y)=0$.
\end{proof}

\begin{lemma}\lbl{L.nc-weak}
Assume $\phi$ and $\psi$ are diffuse Borel 
 Maharam submeasures and  $\phi$  does not have a control measure. Then 
there is a Borel set $A\subseteq 2^\omega\times 2^\omega$ such that 
 $\psi(A_x)=0$ for all $x$ and $\inf_y \phi(A^y)>0$. 
 \end{lemma} 
 
 \begin{proof} 
 Let $A$ be a Borel set
 such that the restriction of $\phi$ to $D^\co$ has a control measure while the restriction 
 of $\phi$ to $D$ is pathological, as given by Lemma~\ref{L.decompose}. 
 By our assumption, $\phi(D)>0$. 
 Again using Lemma~\ref{L.decompose}, find a Borel partition  $2^\omega=B\cup C$ so that 
 $\psi$ has a control measure on $B$ and is pathological on $C$. 
 By Lemma~\ref{L.Chr} there is Borel $E\subseteq D\times B$ such that 
 $\phi(E^y)=\phi(D)$ for all $y\in B$ and $\psi(E_x)=0$ for all $x$. 
 By Lemma~\ref{L.non-Fubini} there is a Borel $F\subseteq D\times C$ 
 such that $\phi(F^y)\geq \frac 13\phi(D)$ for all $y\in C$ 
 and $\psi(F_x)=0$ for all~$x$. Then $A=E\cup F$ is as required. 
 \end{proof} 
 
\section{Non-commutativity}

Given $\gs$-ideals $I$ and $J$ on the real line, let $I\perp J$ be the statement that there is a
Borel subset $B$ of the plane such that all of its vertical sections are in the 
ideal $J$ and all of the
horizontal sections of the complement are in the ideal $I$.
Thus $I\perp J$ means that the Fubini theorem 
between $I$ and $J$ fails in a particularly violent manner.
If the $\gs$-ideals $I$ and $J$ are definable and c.c.c., $I\perp J$ is 
easily seen to be equivalent to both $P_I\Vdash 2^\gw\cap V\in\dot J$
and $P_J\Vdash 2^\gw\cap V\in\dot I$ \cite[5.4.8]{z:book}.

Results of this section do not any require large cardinals if the ideals
are assumed to be analytic on $G_\delta$. This is because in this case both the compatibility 
and the incompatibility relations of $P_I$ are analytic, and therefore the result of \cite{Sh:480} 
applies. Let us recall and  prove Corollary~\ref{C.random.0}.

\begin{corollary}[LC]
\label{C.random}
Suppose  $I$ is a  definable c.c.c. $\gs$-ideal on $2^\gw$. 
Then exactly one of the following holds: 
\begin{enumerate}
\item There is a Borel set $B\subset \R\times\R$
 with all  vertical sections in $I$ and all horizontal sections of 
 full Lebesgue measure. 
\item There is a condition $p\in P_I$ such that $P_I$ below $p$ 
is isomorphic to the random forcing.
\end{enumerate}
\end{corollary}

\begin{proof}
By Fubini's theorem, two clauses exclude each other. 
Assume that $P_I$ is not isomorphic to the random algebra below any positive set $B$. 
By Theorem~\ref{maintheorem} and Lemma~\ref{L.Chr}, we may assume 
   $P_I$ is not bounding, so  by a result of Shelah (\cite{Sh:480})
it adds a Cohen real. 
Let $f\colon 2^\omega\to 2^\omega$ be a Borel function such that the preimages of meager sets
are in $I$. 
Fix a Borel set $B\subseteq 2^\omega\times2^\omega$ whose all 
vertical sections are null and whose complements of horizontal sections are meager. 
Then the set $D=\{(x,y)| (f(x),y)\in A\}$ witnesses that $I\perp\lnull$. 
\end{proof} 

We do not know whether  $\Null(\phi)\perp \Null(\psi)$ whenever
$\phi$ and $\psi$ are Maharam submeasures at least one of which is pathological. 

 Shelah \cite{Sh:630}
defined the notion of commutation for definable c.c.c. $\gs$-ideals
$I,J$: they commute if for all reals $r,s$ in all generic extensions of $V$, the statement
``$r$ is $V[s]$-generic for $P_I$ and $s$ is $V$-generic for $P_J$'' is 
equivalent to ``$s$ is $V[r]$-generic for $P_J$
and $r$ is $V$-generic for $P_I$.'' (Note that in this situation
$r$ is automatically $V$-generic for $P_I$, since it avoids all sets in $I$ coded in $V$.)
Shelah proved that the only ideal commuting with $\meager$ is $\meager$ itself. 
Corollary~\ref{C.random} can be formulated by saying 
 that the only ideal commuting with $\lnull$ is $\lnull$
itself. In \cite{Sh:630}, Problem 11.5, Shelah asked whether the only 
Suslin forcings that commute with themselves
are cohen and random. (A forcing notion $\bbP$ is \emph{suslin} if 
its underlying set is $\R$ and both $\leq_{\bbP}$ and $\perp_{\bbP}$ are 
analytic subsets of the plane. If $I$ is analytic on $G_\delta$, then $P_I$ is easily Suslin.)
By Theorem~\ref{T.non-commuting}, the answer to this question restricted to 
definable forcings of the form $P_I$ is positive.

Rec\l aw and Zakrzewski (\cite{ReZa:Fubini}) 
say that a pair of ideals $I,J$ has the \emph{Fubini property} if 
for every Borel $B\subseteq 2^\omega\times 2^\omega$ such that 
$\{x|B_x\notin J\}=\emptyset$ we have  $\{y|B^y\notin I\}\in J$. 
They have proved that 
 in a certain restricted class of ccc $\sigma$-ideals
of Borel sets $(\meager,\meager)$ and $(\lnull,\lnull)$ are the only pairs  that have the 
Fubini property. They also found a consistent example (using a large cardinal assumption) 
of another ccc  $\sigma$-ideal $I$ such that $(I,I)$ has the  
Fubini property, and asked whether there are other `natural' 
examples of pairs of ccc  ideals with Fubini property. 
Theorem~\ref{T.non-commuting} gives a negative answer to their question restricted to the 
 class of definable ideals.

 In order to give  unified treatment of ideals $\meager$ and $\lnull$ and the corresponding
 forcing notions Cohen and random, in 
 \cite[Definition~1.26]{Ku:Random-Cohen} Kunen introduced the class of `reasonable' ideals. 
Among other properties, every reasonable ideal is 
a \emph{Fubini ideal} 
(\cite[Definition 1.3]{Ku:Random-Cohen}) and  this implies that $(I,I)$ has the Fubini property. 
Therefore by Theorem~\ref{T.non-commuting}, $\meager$ and $\lnull$ are the only reasonable
ideals that are analytic on $G_\delta$. 
The definition of reasonable also involves being absolute 
 (\cite[Definition~1.20]{Ku:Random-Cohen}) and under large cardinals 
 every absolute set of reals belongs to $L(\R)$
by \cite[Theorem 3.2]{z:proper}. 
Therefore  large cardinals imply that $\meager$ and $\null$ are 
are  the only reasonable
ideals. 

By \cite{RoSh:628}  if the assumption that  $I$ is a Fubini ideal is 
dropped from the definition of a reasonable ideal
then there are many  ideals satisfying the weaker notion. 

\begin{lemma} \lbl{L.nc} Suppose $I$ and $J$ are definable c.c.c. $\gs$-ideals on $2^\gw$. 
Then the following are equivalent. 
\begin{enumerate}
\item\lbl{L.nc.1} $P_I$ and $P_J$ commute. 
\item\lbl{L.nc.2}  If $B\subseteq 2^\omega\times 2^\omega$ is Borel then 
$\{x| B_x\notin J\}\notin I $ implies $\{y| B^y\notin I\}\neq \emptyset$. 
\item \lbl{L.nc.3} Pair $J,I$ has the Fubini property.
\end{enumerate}
\end{lemma} 

\begin{proof} 
Assume \eqref{L.nc.2} fails and fix a Borel $B$ such that 
$\{x| B_x\notin J\}\notin I$ and  $C=\{y| B^y\notin I\}\in J$. 
Let  $A=B\setminus 2^\omega\times C$, and note that 
$\{x| A_x\notin J\}\notin I$ and  $\{y| A^y\notin I\}=\emptyset$. 
Let  $x$ be $V$-generic for $P_I$ so that  that $A_x\notin J$
and let $y\in A_x$ be  $V[x]$-generic for $P_J$. 
Since $A^y\in I$ and $x\in A^y$, $x$ is not $P_I$-generic over $V[y]$. 

Now assume \eqref{L.nc.1} fails, and fix a countable transitive model $M$ 
of a large enough fragment of ZFC containing definitions of $I$ and $J$. 
Since $\{x| x$ is $M$-generic for $P_I\}$ is equal to the complement of the 
union of all Borel sets coded in $M$ that belong to $I$, it is Borel. 
Similarly, the  set 
$$
A_{IJ}=\{(x,y)| x\text{ is $M$-generic for $P_I$ and $y$ is $M[x]$-generic for $P_J$}\}
$$
is Borel, and $B=A_{IJ}\setminus \{(x,y)| (y,x)\in A_{JI}\}$ is a 
Borel set consisting of all pairs $(x,y)$ that fail the commutativity condition. 
This set is nonempty by our assumption, and it satisfies \eqref{L.nc.2}.

To see that \eqref{L.nc.2} and \eqref{L.nc.3} are equivalent, 
take the contrapositive of \eqref{L.nc.2}. 
\end{proof} 

\begin{theorem}[LC] \lbl{T.non-commuting}
Suppose  $I,J$ are definable c.c.c. $\gs$-ideals on $2^\gw$. Then one of the following holds:
\begin{enumerate}
\item Both $P_I$ and $P_J$ are isomorphic to the Cohen algebra. 
\item Both $P_I$ and $P_J$ are isomorphic to the Lebesgue measure algebra. 
\item $P_I$ and $P_J$ do not commute. 
\end{enumerate}
In particular, if $P_I$ of this kind commutes with itself, then it is either Cohen or random. 
\end{theorem}

The proof of Theorem~\ref{T.non-commuting}
breaks into several cases according to whether the posets $P_I, P_J$ are bounding or not,
with wildly different arguments in each case.

\begin{lemma}[LC]
\label{unboundedlemma}
 Suppose  $I,J$ are definable c.c.c. $\gs$-ideals on $2^\gw$ such that both forcings $P_I$ and $P_J$
add an unbounded real. Exactly one of the following holds:
\begin{itemize}
\item there are Borel $I$-positive set $B$ and a Borel $J$-positive 
set $C$ such that both $P_I$ below $B$
and $P_J$ below $C$ are isomorphic to the Cohen algebra
\item $I\perp J$
\end{itemize}
\end{lemma}

\begin{proof}
There is nothing really new here. Clearly the first item implies the failure of $I\perp J$. 
On the other hand,
suppose that the first item fails. 
Then one of the partial orders, $P_I$ say, is not isomorphic to the Cohen algebra 
below any condition.
By \cite{Sh:630} 9.16 or
\cite{z:products} 6.6, $P_I\Vdash 2^\gw\cap V$ is meager, so $I\perp\meager$ and there 
is a Borel set
$E\subset 2^\gw\times 2^\gw$ such that its vertical sections are meager and the horizontal 
sections of its complement
are $I$-small. By \cite{Sh:480}, 1.14, $P_J$ adds a Cohen real over $V$ and so there is 
a Borel function $f:2^\gw\to 2^\gw$ such that
preimages of meager sets are $J$-small. It is not difficult to verify that the 
Borel set $D\subset2^\omega\times2^\omega$
defined by $\langle x, y\rangle\in D$ if and only 
if $\langle x, f(y)\rangle\in E$ witnesses $I\perp J$.
The lemma follows.
\end{proof}

\begin{lemma}[LC] \lbl{L.wd} Suppose  $I,J$ are definable c.c.c. $\gs$-ideals 
such that both forcings $P_I$ and $P_J$ are bounding.
If   $P_I$ is not equivalent to random,
then there is a Borel  $B\subseteq 2^\omega\times 2^\omega$ such that 
$B_x\in J$ for all $x$ and $B^y\notin I$ for all $y$. 
\end{lemma}

\begin{proof} By Theorem~\ref{maintheorem}, both $I$ and $J$ are null ideals 
for  some Maharam submeasures $\phi$ and $\psi$, respectively. 
By Lemma~\ref{L.path},  $\phi$ does not have a control measure. 
Therefore we are in the situation of    Lemma~\ref{L.nc-weak}.
\end{proof}

\begin{lemma}
[LC]\lbl{L3} Suppose that $I,J$ are definable c.c.c. $\gs$-ideals on
 $2^\gw$ such that $P_I$ is bounding while $P_J$ adds
an unbounded real. Then $I\perp J$.
\end{lemma}

\begin{proof} By Theorem~\ref{maintheorem}, there is a Maharam submeasure $\phi$ such that 
$I$ is the null ideal for $\phi$. 
We will first prove that $I\perp \meager$. 
For $s\in 2^n$ let $[s]=\{x\in 2^\omega| x\restriction n=s\}$. 

\begin{claim} If $\phi$ is a Maharam submeasure on the Borel algebra of 
$2^\omega$, then for every $\e>0$ there is $m_\e\in \N$ such that $\phi([s])\leq \e$
for every $s\in 2^{m_\e}$. 
\end{claim} 

\begin{proof} Assume not, and find $s_m\in 2^m$ such that $\phi([s_m])\geq \e$ for all $m$.
Ramsey's theorem gives us two possibilities.
 
Either there is an infinite set $B\subseteq \omega$ such that $[s_m]$ ($m\in B$) 
are pairwise disjoint. In this case the open sets $U_n=\bigcup\{ [s_m]| m\geq n, n\in B\}$ 
have all submeasure at least $\e$ and they are
decreasing with empty intersection. 
Since $\phi$ is a Maharam submeasure, this is impossible.

Or there is an infinite set $D$
such that $[s_m]$ ($m\in D$) form a decreasing chain. The 
intersection $\bigcap_{m\in D}[s_m]$ is a singleton, $\{x\}$,
and again by the continuity of the submeasure, $\phi(\{x\})\geq\e$. 
Thus $\{x\}\notin I$, contradiction.
\end{proof} 

Let $f(n)=m_{2^{-n}}$ as given by the previous Claim. 
Interpret the Cohen forcing as adding a function $g\in\prod_n 2^{f(n)}$ with finite conditions.
Let $D_m=\bigcup_{n>m}[g(n)]$. It is not difficult to see that $V\cap 2^\gw\subset D_m$ 
for every number
$m\in\gw$ and the submeasures $\phi(D_m)$ converge to zero. Therefore $\bigcap_m D_m$ 
is a submeasure zero
set containing all the ground model reals.

Finally, to show that $I\perp J$ note that by a result of Shelah \cite{Sh:480} 
the poset $P_J$ adds a Cohen real.
The argument is concluded in a manner similar to Lemma~\ref{unboundedlemma}.
\end{proof} 

\begin{proof}[Proof of Theorem~\ref{T.non-commuting}] 
Let $I,J$ be definable c.c.c. $\gs$-ideals,
and suppose that the first two alternatives in the Theorem fail. Use the c.c.c.
to find partitions $2^\gw=B_0\cup B_1$ and $2^\gw=C_0\cup C_1$ into Borel sets such 
that $P_I$ below $B_0$
and $P_J$ below $C_0$ are bounding forcings while the posets $P_I$ below $B_1$ and 
$P_J$ below $C_1$ add
an unbounded real. Pick $i,j$ such that $B_i\notin I$ and $C_j\notin J$. 
If $i=j$ we may assure that if $P_I$ is $\meager$ ($\lnull$, respectively) below $B_i$ 
then  $P_J$ is not $\meager$ ($\lnull$, respectively) below $C_j$. 
In either case, by one of lemmas \ref{unboundedlemma}, \ref{L.wd} or \ref{L3} 
we are in the situation of Lemma~\ref{L.nc}. 
\end{proof} 

\section{Concluding remarks}

Another corollary of Theorem~\ref{maintheorem}  precisely determines 
the extent of ccc-ness of a weakly distributive definable forcing $P_I$. 
 Recall that a subset $F$ of a  poset $\bbP$ is $n$-linked 
if every $n$-element subset of $F$ has a lower bound, 
and that $\bbP$ is \emph{$\sigma$-$n$-linked} if it can be covered 
by countably many $n$-linked sets. 
An $F\subseteq \bbP$ is centered if every finite subset of $F$ has a 
lower bound, and $\bbP$ is $\sigma$-centered if it can be covered by countably
many centered subsets. It is well-known that all these chain conditions are different. Also, 
by a result of Todorcevic (\cite{To:Examples}, see also \cite[3.6.C]{BarJu:Book}), 
there is a  Borel ccc 
poset that is not $\sigma$-2-linked.

\begin{corollary} [LC] If $I$ is a  $\sigma$-ideal of Borel sets and 
$P_I$ is weakly distributive, then the following hold. 
\begin{enumerate}
\item\label{ccc.1} $P_I$ is not $\sigma$-centered. 
\item\label{ccc.2} If $I$ is moreover definable, 
then  $P_I$ is ccc if and only if it is $\sigma$-$n$-linked for all $n$. 
\end{enumerate}
\end{corollary}

\begin{proof} 
Assume $\B$ is $\sigma$-centered and fix centered sets $X_n$ maximal under the inclusion 
whose union covers $\B$.  
Since by Fact~\ref{zfact} 
every positive set has a compact subset the intersection of each $X_n$ is a singleton. 
This 
implies that a co-countable set belongs to $I$, a contradiction. 

Only clause~\ref{ccc.2} requires the definability and large cardinal assumptions. 
By Theorem~\ref{maintheorem} it suffices to prove a well-known fact that 
if $\phi$ is a Maharam submeasure on Borel algebra of $2^\omega$ and $\Null(\phi)$ contains
all countable sets, 
 then 
the quotient algebra is $\sigma$-$n$-linked for all $n$
(this is \cite[Exercise 393Y(a)]{Fr:MT3}). 
 Recall first that it 
 is completely generated by its countable subalgebra $\B_0$ given by the name 
 for the $P_I$-generic real. 
Now consider the metric on $\B$ defined by $\rho(A,B)=\phi(A\Delta B)$. 
It is not difficult to check that $(\B,\rho)$ is a complete metric space, 
and as an easy consequence of 
 \cite[393B (c)]{Fr:MT3}, it is isomorphic to the  completion of $(\B_0,\rho)$, 
and in particular separable. 
For $A\in \B_0$ the set 
$$
F_A=\{C|\rho(A,C)<\rho(0_\B,A)/n\}. 
$$
is $n$-linked, and $\bigcup_{A\in \B_0} F_A$ covers $\B$. 
\end{proof}
 
We conclude with a question asked by Solecki (personal communication). 

\begin{question} Are the following equivalent for 
every  c.c.c. $\gs$-ideal $I$ on Borel subsets of $2^\gw$ that is 
analytic on $G_\delta$?
\begin{enumerate}
\item Compact sets are dense in $P_I$ and $I$ is ccc. 
\item $I$ is the null ideal of some Maharam submeasure. 
\end{enumerate}
\end{question} 

If the answer is positive, this would strengthen Theorem~\ref{maintheorem} and nicely complement
a result of \cite{KeSol:Approximating} where it was proved that every ccc $\sigma$-ideal 
$\sigma$-generated by compact sets is Borel-isomorphic to  $\meager$.

\def\germ{\frak} \def\scr{\cal} \ifx\documentclass\undefinedcs
  \def\bf{\fam\bffam\tenbf}\def\rm{\fam0\tenrm}\fi % f**k-amstex!
  \def\defaultdefine#1#2{\expandafter\ifx\csname#1\endcsname\relax
  \expandafter\def\csname#1\endcsname{#2}\fi} \defaultdefine{Bbb}{\bf}
  \defaultdefine{frak}{\bf} \defaultdefine{mathbb}{\bf}
  \defaultdefine{mathcal}{\cal}
  \defaultdefine{beth}{BETH}\defaultdefine{cal}{\bf} \def\bbfI{{\Bbb I}}
  \def\mbox{\hbox} \def\text{\hbox} \def\om{\omega} \def\Cal#1{{\bf #1}}
  \def\pcf{pcf} \defaultdefine{cf}{cf} \defaultdefine{reals}{{\Bbb R}}
  \defaultdefine{real}{{\Bbb R}} \def\restriction{{|}} \def\club{CLUB}
  \def\w{\omega} \def\exist{\exists} \def\se{{\germ se}} \def\bb{{\bf b}}
  \def\equivalence{\equiv} \let\lt< \let\gt> \def\cite#1{[#1]}


\begin{thebibliography}{10}

\bibitem{BJP}
Balcar. B., T.~Jech, and T.~Paz{\'a}k.
\newblock Complete ccc boolean algebras, the order sequential topology, and a
  problem of von {N}eumann.
\newblock preprint, 2003.

\bibitem{BarJu:Book}
T.~Bartoszynski and H.~Judah.
\newblock {\em Set theory: on the structure of the real line}.
\newblock A.K. Peters, 1995.

\bibitem{Chr:Fubini}
J.P.R. Christensen.
\newblock Some results with relation to the control measure problem.
\newblock In R.M. Aron and S.~Dineen, editors, {\em Vector space measures and
  applications {II}}, volume 645 of {\em Lecture Notes in Mathematics}, pages
  27--34. Springer, 1978.

\bibitem{Fr:MT3}
D.H. Fremlin.
\newblock {\em Measure Theory}, volume~3.
\newblock Torres--Fremlin, 2002.

\bibitem{HeChr:Example}
W.~Herer and J.P.R. Christensen.
\newblock On the existence of pathological submeasures and the construction of
  exotic topological groups.
\newblock {\em Mathematische Annalen}, 213:203--210, 1975.

\bibitem{z:proper}
{Itay Neeman} and Jind{\v r}ich Zapletal.
\newblock Proper forcing and absoluteness in l(r).
\newblock {\em Commentationes Mathematicae Universitatis Carolinae},
  39:281--301, 1998.

\bibitem{jech:set}
Thomas Jech.
\newblock {\em Set Theory}.
\newblock Academic Press, San Diego, 1978.

\bibitem{Ka:Maharam}
N.J. Kalton.
\newblock The {M}aharam problem.
\newblock In G.~Choquet et~al., editors, {\em S\'eminaire Initiation \`a \'l
  Analyse, 28e An\'nee}, volume~18, pages 1--13. 1988/89.

\bibitem{KaRo:Uniformly}
N.J. Kalton and J.W. Roberts.
\newblock Uniformly exhaustive submeasures and nearly additive set functions.
\newblock {\em Transactions of the American Mathematical Society},
  278:803--816, 1983.

\bibitem{Kana:Book}
A.~Kanamori.
\newblock {\em The higher infinite: large cardinals in set theory from their
  beginnings}.
\newblock Perspectives in Mathematical Logic. Springer--Verlag,
  Berlin--Heidelberg--New York, 1995.

\bibitem{Ke:Classical}
A.S. Kechris.
\newblock {\em Classical descriptive set theory}, volume 156 of {\em Graduate
  texts in mathematics}.
\newblock Springer, 1995.

\bibitem{KeSol:Approximating}
A.S. Kechris and S.~Solecki.
\newblock Approximating analytic by {B}orel sets and definable chain
  conditions.
\newblock {\em Israel J. Math.}, 89.

\bibitem{Ku:Random-Cohen}
K.~Kunen.
\newblock Random and cohen reals.
\newblock In K.~Kunen and J.~Vaughan, editors, {\em Handbook of Set-Theoretic
  Topology}. North-Holland, 1984.

\bibitem{Lo:Progres}
A.~Louveau.
\newblock Progres recents sur le probleme de {M}aharam d'apres {N.J.} {K}alton
  et {J.W.} {R}oberts.
\newblock In G.~Choquet et~al., editor, {\em S\'eminaire Initiation \`a \'l
  Analyse, 28e An\'nee}, volume~20, pages 01--08. 1983/84.

\bibitem{Mah:Algebraic}
D.~Maharam.
\newblock An algebraic characterization of measure algebras.
\newblock {\em Annals of mathematics}, 48:154--167, 1947.

\bibitem{Mah:vonNeumann}
D.~Maharam.
\newblock Problem 167.
\newblock In D.~Mauldin, editor, {\em The {S}cottish {B}ook}, pages 240--243.
  Birk\-h\"au\-ser, Boston, 1981.

\bibitem{ReZa:Fubini}
Ireneusz Rec{\l}aw and Piotr Zakrzewski.
\newblock Fubini properties of ideals.
\newblock {\em Real Anal. Exchange}, 25(2):565--578, 1999/00.

\bibitem{RoSh:628}
Andrzej Roslanowski and Saharon Shelah.
\newblock Norms on possibilities ii: more ccc ideals on $2^{\textstyle\omega}$.
\newblock {\em Journal of Applied Analysis}, 3:103--127, 1997.

\bibitem{Sh:630}
Saharon Shelah.
\newblock Properness without elementaricity.
\newblock {\em Journal of Applied Analysis}, accepted.

\bibitem{Sh:480}
Saharon Shelah.
\newblock How special are cohen and random forcings i.e. boolean algebras of
  the family of subsets of reals modulo meagre or null.
\newblock {\em Israel Journal of Mathematics}, 88:159--174, 1994.

\bibitem{Sol:AnalyticII}
S.~Solecki.
\newblock Analytic ideals and their applications.
\newblock {\em Annals of Pure and Applied Logic}, 99:51--72, 1999.

\bibitem{To:Definable}
S.~Todorcevic.
\newblock Definable ideals and gaps in their quotients.
\newblock In {C.A.} DiPrisco et~al., editors, {\em Set Theory: {T}echniques and
  Applications}, pages 213--226. Kluwer Academic Press, 1997.

\bibitem{To:Dichotomy}
S.~Todorcevic.
\newblock A dichotomy for {P}-ideals of countable sets.
\newblock {\em Fundamenta Mathematicae}, 166:251--267, 2000.

\bibitem{To:Examples}
Stevo Todor{\v{c}}evi{\'c}.
\newblock Two examples of {B}orel partially ordered sets with the countable
  chain condition.
\newblock {\em Proc. Amer. Math. Soc.}, 112(4):1125--1128, 1991.

\bibitem{z:products}
J.~Zapletal.
\newblock Proper forcing and rectangular {R}amsey theorems.
\newblock preprint.

\bibitem{z:book}
Jind{\v r}ich Zapletal.
\newblock {\em Descriptive Set Theory and Definable Forcing}.
\newblock Memoirs of American Mathematical Society. AMS, Providence, 2004.

\end{thebibliography}
\end{document}